\documentclass[11pt]{article}
\usepackage{mathrsfs}
\usepackage{amssymb}
\usepackage{amsfonts}
\usepackage{color,xcolor}
\usepackage{graphicx}
\usepackage{geometry}
\usepackage{manfnt}
\usepackage[pagewise]{lineno}
\usepackage[colorlinks,citecolor=blue,urlcolor=blue]{hyperref}
\newtheorem{theorem}{Theorem}[section]

\newtheorem{lemma}{Lemma}[section]

\newtheorem{remark}{Remark}[section]
\newtheorem{example}{Example}[section]

\newcommand{\beq}{\begin{equation}}
\newcommand{\eeq}{\end{equation}}
\newcommand{\beqn}{\begin{eqnarray}}
\newcommand{\eeqn}{\end{eqnarray}}

\def\[{{\Big[}}\def\]{{\Big]}}\def\<{{\langle}}\def\>{{\rangle}}\def\({{\Big(}}
\def\){{\Big)}}

\def\min{{\mathord{{\rm min}}}}
\def\={&\!\!=\!\!&}

\def\geq{\geqslant}\def\leq{\leqslant}

\begin{document}
\title{\bf
On mountain pass theorem and its application to periodic solutions of some nonlinear discrete systems\thanks{E-mail: ding2016liang@126.com (L. Ding), weijinlong@zuel.edu.cn (J. Wei), zhangshiqing@msn.com (S. Zhang).}
\,\,\thanks{
The project is supported by the National Natural Science Foundation of China (grant no. 1150157 and 11671278).}}

\author{Liang Ding$^a$, Jinlong Wei$^b$, Shiqing Zhang$^a$
\date{{\small \it $^a$ School of Mathematics, Sichuan University, Chengdu 610064, China}\\
{\small \it $^b$ School of Statistics and Mathematics, Zhongnan University }\\ {\small \it of Economics and Law, Wuhan 430073, China}}}

 \maketitle
\noindent{\hrulefill} \vskip1mm\noindent
 {\bf Abstract}
 \vskip1mm\noindent
 We obtain a new quantitative deformation lemma, and then gain a new mountain pass theorem. More precisely, the new mountain pass theorem is independent of the functional value on the boundary of the mountain, which improves the well known results (\cite{AR,PS1,PS2,Qi,Wil}). Moreover, by our new mountain pass theorem, new existence of  nontrivial periodic solutions for some nonlinear second-order discrete systems is obtained, which greatly improves the result in \cite{Z04}.
  \vskip2mm\noindent
{\bf Mathematics Subject Classification.} 49J35, 34K13, 34G20
\vskip2mm\noindent
{\bf Keywords.} Critical points, Quantitative deformation lemma,
Mountain pass theorem, Nontrivial periodic solutions, Nonlinear discrete systems
\vskip1mm\noindent{\hrulefill}

\section{Introduction and main results} \label{sec1}
\setcounter{equation}{0}
It is well known that the classical mountain pass theorem of Ambrosetti-Rabinowitz \cite{AR} has been proved to be a powerful tool in applications to many areas of analysis, and quantitative deformation lemma is used to be a very excellent method to derive different kinds of minimax theorems including the mountain pass theorem, we refer the authors to see \cite{Cha,Gho,Sch,Wil}. Firstly, we recall the famous mountain pass theorem given by Ambrosetti and Rabinowitz \cite{AR}:
\begin{theorem} [\cite{AR}] \label{the1.1}
Let $X$ be a Banach space, $\varphi\in C^{1}(X, \mathbb{R})$, suppose there exist  $e\in X$ and two real numbers $\alpha>0$ and $r>0$ such that
  $\|e\|>r$ and   \begin{itemize}
  \item[$(i)$] $c_0\geq \alpha>0$ where $c_0:=\inf_{\|u\|=r}\varphi(u)$, and $\varphi(u)>0$ in $\{u\in X, \  |\|u\|<r\} \setminus\{0\}$;
 \item[$(ii)$] $\varphi(0)=\varphi(e)=0$;
\item[$(iii)$]If $(u_{n}) \subset X$ with $0<\varphi(u_{n}), \,\,\varphi(u_{n})$ bounded above, and
 $\varphi^{\prime}(u_{n})\rightarrow 0$, then $(u_{n})$ possesses a convergent subsequence.
\end{itemize}
Then, $c:=\inf_{\gamma\in\Gamma}\max_{t\in[0,
      1]}\varphi\big(\gamma(t)\big)$,
where
\begin{eqnarray*}
\Gamma:=\{\gamma\in C\big([0, 1], X\big):\gamma(0)=0, \,\,\gamma(1)=e\},
\end{eqnarray*}
is a critical value of $\varphi$.
\end{theorem}

Set $c_1:=\max\{\varphi(0), \,\,\varphi(e)\}$, we find that $c_0>c_1$ in Theorem \ref{the1.1} and since then, there are many variant generalizations on this case for mountain pass theorem \cite{AR,BCN,Cha,Wil}. One of the elegant works is founded by Willem \cite{Wil}. To be make the result more clear, we outline the quantitative deformation lemma first
\begin{lemma} [\cite{Wil}] \label{lem1.1} (Quantitative deformation lemma)
 Let $X$ be a Hilbert space, $\varphi\in C^{2}(X, \mathbb{R})$, $c\in \mathbb{R}$, $\varepsilon>0$. Assume that
\begin{eqnarray*}
\|\varphi^{\prime}(u)\|\geq2\varepsilon, \quad \forall \ u\in\varphi^{-1}([c-2\varepsilon,
  c+2\varepsilon]).
\end{eqnarray*}
Then there exists $\eta\in$ $C(X, X)$, such that
\begin{itemize}
\item[$(i)$] $\eta(u)=u$, $\forall\ u\notin \varphi^{-1}\big([c-2\varepsilon,
c+2\varepsilon]\big)$;
\item[$(ii)$]
$\eta(\varphi^{c+\varepsilon})\subset\varphi^{c-\varepsilon}$,
where $\varphi^{c-\varepsilon}:=\varphi^{-1}\big((-\infty, c-\varepsilon]\big)$.
\end{itemize}
\end{lemma}
If one lets $c_0$ be stated in Theorem \ref{the1.1}, then Willem's result can be described as the following:
\begin{theorem} [\cite{Wil}] \label{the1.2} (Mountain pass type theorem)
  Let $X$ be a Hilbert space, $\varphi\in C^{2}(X, \mathbb{R})$,
suppose
\begin{itemize}
\item[$(i)$]$\varphi$ satisfies the $(P.S.)$ condition ($\varphi$ is said to satisfies $(P.S.)$ condition, if any sequence $\{u^{(k)}\}\subset X$ satisfying $\varphi(u^{(k)})\rightarrow c$ and $\varphi^{\prime}(u^{(k)})\rightarrow 0$ as $k\rightarrow+\infty$ with any constant $c$, implies a convergent subsequence);
\item[$(ii)$]there exist $e\in X$ and $r>0$ be such that $\|e\|>r$ and $c_0>\varphi(0)\geq\varphi(e)$.
\end{itemize}
Then, $c$ is a critical value of $\varphi$.
\end{theorem}

For the case of $c_0\geq c_1$ in mountain pass type theorems, we refer the authors to see \cite{PS1,PS2,Qi}. Specially, we introduce a mountain pass type theorem in \cite{PS2} and the extension mountain pass type theorem in \cite{Qi} as follows:
\begin{theorem} [\cite{PS2}] \label{the1.3} (Mountain pass type theorem)
  Let $X$ be a Banach space and $X$ has finite dimension, $\varphi\in C^{1}(X, \mathbb{R})$. Suppose there exist  $e\in X$ and two real numbers $a$ and $r>0$ such that
  $\|e\|>r$ and
 \begin{itemize}
  \item[$(i)$] $c_0\geq a$;
 \item[$(ii)$] $\varphi(0)\leq a$, $\varphi(e)\leq a$;
 \item[$(iii)$] any sequence $(u_{n})$ in $X$ such that $\varphi(u_{n})\rightarrow limit \geq a$, and $\varphi^{\prime}(u_{n})\rightarrow 0$
 possesses a convergent subsequence.
\end{itemize}
Then, $c$ is a critical value of $\varphi$.
\end{theorem}
\begin{remark}
For infinite-dimensional case, Pucci and Serrin in \cite{PS2} gained that $c$ is a critical value of $\varphi$ when the assumption $(i)$ in Theorem \ref{the1.3} was suitably strengthened, more precisely, their conditions depending on $c_0$ and the neighbor of $\{\varphi(u), \ \|u\|=r\}$.
 \end{remark}
\begin{theorem} [\cite{Qi}] \label{the1.4}(Extension mountain pass theorem)
Let $X$ be a real Hilbert space, $\varphi\in C^{1}(X, \mathbb{R})$ satisfying the $(P.\,S.)$ condition, $e\in X$ and $r>0$ be such that
$\|e\|>r$.  If $c_0\geq c_1$, $c$ is a critical value of $\varphi$.
\end{theorem}

Evidently, all above Theorems are based upon $c_0$. Then, an interesting question is raised: can we obtain a new mountain pass theorem which is independent of $c_0$ ? In this paper, we give a positive answer and the new mountain pass theorem is given by:
\begin{theorem}\label{the1.5}(New mountain pass theorem) Let $X$ be a Hilbert space,
$\varphi\in C^{2}(X, \mathbb{R})$, $e, \,e_1\in X$,
 $r>0$ be such that $0<\|e_1\|<r$ and $\|e\|>r$, and
$\varphi(0)<\varphi(e)=\varphi(e_1)$.
Then, for each $\varepsilon>0$, there exists $\hat{u}\in X$ such that
\begin{itemize}
\item[$(i)$] $\hat{c}-2\varepsilon\leq\varphi(\hat{u})\leq \hat{c}+2\varepsilon$;
\item[$(ii)$] $\|\varphi^{\prime}(\hat{u})\|<2\varepsilon$,
\end{itemize}
where $\hat{c}:=\inf_{\gamma\in\hat{\Gamma}}\max_{t\in[0, 1]}\varphi\big(\gamma(t)\big)$ and
\begin{eqnarray*}
\hat{\Gamma}:=\{\gamma\in C\big([0, 1], X\big):\gamma(0)=0, \gamma(\frac{1}{2})=e_1,\gamma(1)=e\}.
\end{eqnarray*}
\end{theorem}
\begin{remark}
The new mountain pass theorem is independent of $c_0$, and if $\varphi$ satisfies the $(P.S.)$ condition, there exists $\hat{u}\in X$ such that $\varphi(\hat{u})= \hat{c}$ .
\end{remark}

Now, we turn to an application of our new mountain pass theorem to the existence of periodic solutions of discrete systems, which has been  appeared in computer science, economic, neural networks, ecology,
cybernetics, etc and extensively investigated in \cite{Ch,De,D17,GY1,GY2,Ra,X07,X08,Y06,Z04}.

Let $\mathbb{Z}, \,\mathbb{N}, \,\mathbb{R}$ be the set of all integers, natural numbers and real numbers, respectively.
In \cite{GY1}, by critical point theory, Guo and Yu established the existence of periodic solutions to the below discrete difference equations
\begin{eqnarray}\label{1.1}
\Delta^{2}u_{n-1}+f(n,\,u_{n}) = 0, \quad  n \in \mathbb{Z},
 \end{eqnarray}
where $f(n, u_n)=\nabla_{u_n} F(n, u_n)$, $u_n=u(n)\in\mathbb{R}, \ \Delta u_n = u_{n+1}-u_n$,
$\Delta^{2}u_n=\Delta(\Delta u_n)$, $F:
\mathbb{Z}\times\mathbb{R}\rightarrow \mathbb{R}$, $F(n,x)$ is
continuously differentiable in $x$ for every $n\in \mathbb{Z}$ and
$T$-periodic ($0<T\in \mathbb N$) in $n$ for all
$x\in\mathbb{R}$, $\nabla_x F(n,x)$ is the gradient of $F(n,x)$
in $x$. So far as we know, \cite{GY1} is the first paper to study the existence of periodic solutions of system (\ref{1.1}) for superlinear $f(n,u_n)$. For more results when $f(n, u_n)$ is superlinear in the second variable $u_n$, one consults to \cite{D17,X08}.
When $f(n, u_n)$ is sublinear in the second variable $u_n$, we refer the authors to see \cite{GY2} and \cite{X07} and for the case of $f(n, u_n)$ is neither suplinear nor sublinear, we refer to see \cite{Z04}. For more details in this direction, one consults to \cite{Ch,De,Ra,Y06}.
It is remarked that, in \cite{Z04}, under the assumptions described below:
\begin{itemize}
\item [$(A_{1})$]$F\geq0, \,F\in C^{1}(\mathbb{R}\times\mathbb{R}, \,\mathbb{R})$, and for every $(n, \,x)\in\mathbb{Z}\times\mathbb{R}$, there is a positive integer $M\geq3$ such that $F(n+M,x)=F(n,x)$;
\item [$(A_{2})$]there exist constants $\delta>0$, $\alpha\in(0,1-\cos\frac{2\pi}{M})$ such that
\begin{eqnarray*}
F(n,\,x)\leq\alpha x^{2} \ \  \mbox{for} \ \ n\in \mathbb{N},\
x\in\mathbb{R} \ \mbox{and} \ |x|\leq\delta;
\end{eqnarray*}
\item [$(A_{3})$]
there exist constants $w_1>0$, $w_2>0$ and
$w_3\in(2,\,+\infty)$ when $M$ is even or $w_3\in(1+\cos\frac{\pi}{M},\,+\infty)$ when $M$ is odd, such that
\begin{eqnarray*}
F(n, \,x) \geq w_3x^{2}-w_2 \ \ \mbox{for} \ n\in\mathbb{N}, \ |x|\geq w_1.
\end{eqnarray*}
\end{itemize}
by using linking theorem \cite{Ra}, Zhou, Yu and Guo derived the existence of
nontrivial $M$-periodic solutions for system (\ref{1.1}), and they give an example:

\begin{example}\label{exa1.1}
Take $F(t,\,x)=a(x^{2}/2+\cos x-1)(\phi(t)+K)$ with $x\in\mathbb{R}$, constant $K>0$, where constant $a$ and
positive integer $M\geq3$ satisfy
\begin{eqnarray*}
\cases{a>2, \ \   \mbox{when} \  M \  \mbox{is even}, \cr
a>2(1+\cos\frac{\pi}{M}), \  \mbox{when}  \ M \  \mbox{is odd,}}
\end{eqnarray*}
$\phi(t)\in C^{1}({\mathbb{R},\mathbb{R}})$, and $\phi(t)$ is a $M$-periodic function satisfying $|\phi(t)|<K$.
\end{example}

Obviously, $F(n, \,0)=0, \,\,\forall n\in \mathbb{Z}$, is weaker than condition $(A_{2})$, and the following example is failure to satisfy condition $(A_{2})$, but satisfy $F(n, \,0)=0, \,\,\forall n\in \mathbb{Z}$:
\begin{example}\label{exa1.2}
Take $F(t,\,x)=a(\mu x^{2}+\cos x-1)(\phi(t)+K)$ with $x\in\mathbb{R}$, constants $\mu>1/2$, $K>0$, where constant $a$ and
positive integer $M\geq6$ satisfy
\begin{eqnarray*}
\cases{a>2, \ \   \mbox{when} \  M \  \mbox{is even}, \cr
a>2(1+\cos\frac{\pi}{M}), \  \mbox{when}  \ M \  \mbox{is odd,}}
\end{eqnarray*}
 $\phi(t)\in C^{2}({\mathbb{R},\mathbb{R}})$, and $\phi(t)$ is a $M$-periodic function satisfying $|\phi(t)|<K$.
\end{example}

So, the second interesting question is raised:
for $f(n, u_n)$ is neither superlinear nor sublinear, when condition $(A_{2})$ is replaced by $F(n, \,0)=0, \,\,\forall n\in \mathbb{Z}$,
can we still obtain the existence of nontrivial periodic solutions ?

In this paper, employing our new mountain pass theorem, we obtain new existence of nontrivial periodic solutions for discrete second-order discrete system (\ref{1.1}), and our result is that:
\begin{theorem}\label{the1.6}
Assume that $F\in C^{2}(\mathbb{R}\times\mathbb{R}, \,\mathbb{R})$ and there is a positive integer $M\geq6$ satisfying condition $(A_{3})$ and the following conditions:
\begin{itemize}
\item [$(W_{1})$]$F\geq0$, and for every $(n,x)\in\mathbb{Z}\times\mathbb{R}$, $F(n+M,x)=F(n,x)$;
\item [$(W_{2})$]$F(n,0)=0, \,\,\forall \ n\in \mathbb{Z}$.
\end{itemize}
 Then, system (\ref{1.1}) has  at least one
nontrivial $M$-periodic solutions.
\end{theorem}

\begin{remark}\label{1.3}
Obviously, condition $(W_{2})$ is weaker than condition $(A_{2})$, and it is easy to verify that $F(t,x)$ in Example \ref{exa1.2} satisfies all the conditions of Theorem \ref{the1.6}, but does not satisfy condition $(A_2)$. So, by a different method, i.e., Theorem \ref{the1.5}, we greatly improve the result in \cite{Z04}.
\end{remark}

The paper is organized as follows: Section 2 is devoted to establish a new
quantitative deformation lemma. In Section 3, by using the new quantitative deformation lemma, we derive our new mountain pass theorem (Theorem \ref{the1.5}). In Section 4, as an application of our new mountain pass theorem, we prove Theorem \ref{the1.6}.

\section{New quantitative deformation lemma}
\setcounter{equation}{0}
\begin{lemma}\label{lem2.1} Let $X$ be a Hilbert space, $\varphi\in C^{2}(X, \mathbb{R})$, $h\in \mathbb{R}$,
$\varepsilon>0$. Assume that
\begin{eqnarray*}
\|\varphi^{\prime}(u)\|\geq2\varepsilon, \quad \forall \  u\in\varphi^{-1}([h-2\varepsilon,
h+2\varepsilon]).
\end{eqnarray*}
Then there exists $\eta\in$ $C(X, X)$, such that
\begin{itemize}
\item[$(i)$] $\eta(u)=u$, $\forall \ u\notin \varphi^{-1}\big([h-2\varepsilon,
h+2\varepsilon]\big)\backslash D$, where $D$ is a
subset of $X$ satisfying
$D\subset\varphi^{-1}\big([h-\frac{1}{3}\varepsilon,
h+\frac{1}{3}\varepsilon]\big)$;
\item[$(ii)$] $\eta\big(\varphi^{-1}[h-\varepsilon, h-\frac{1}{2}\varepsilon]\big)\subset
\varphi^{-1}\big([h+\varepsilon, h+\frac{3}{2}\varepsilon]\big)$;
\item[$(iii)$] $\eta\big(\varphi^{-1}[h+\frac{1}{2}\varepsilon, h+\varepsilon]\big)\subset\varphi^{-1}\big([h-\frac{3}{2}\varepsilon,
h-\varepsilon]\big)$.
\end{itemize}
\end{lemma}
\textbf{Proof.} Let us define
\begin{eqnarray*}
&A&:=\varphi^{-1}\big([h-2\varepsilon, h+2\varepsilon]\big)\backslash D,\ \ B:=\varphi^{-1}\big([h-\varepsilon, h-\frac{1}{2}\varepsilon]\big),\\
&C&:=\varphi^{-1}\big([h+\frac{1}{2}\varepsilon, c+\varepsilon]\big),\\
&\psi(u)&:=\frac{[dist(u, C)-dist(u, B)]dist(u, X\backslash
A)}{[dist(u, C)+dist(u, B)]dist(u, X\backslash A)+dist(u, B)dist(u,
C)},
\end{eqnarray*}
so that $\psi$ is locally Lipschitz continuous, $\psi=1$ on $B$,
$\psi=-1$ on $C$ and
$\psi=0$ on $X\backslash$ $A$.
 \vskip1mm\par
Let us also define the locally Lipschitz continuous vector field
\begin{eqnarray}\label{2.1}
f(u):=\cases{\psi(u)\|\nabla\varphi(u)\|^{-2}\nabla\varphi(u), \quad u\in A,\cr \quad 0, \quad \ \  u\in X\backslash A.}
\end{eqnarray}
It is clear that $\|f(u)\|\leq (2\varepsilon)^{-1}$ on $X$. For each
$u\in$ $X$, the Cauchy problem
\begin{eqnarray*}
\frac{d}{dt}\sigma(t, u)=f\big(\sigma(t, u)\big),\ \
\sigma(0, u)=u,
\end{eqnarray*}
has a unique solution $\sigma(\cdot, u)$ defined on $\mathbb{R}$.
Moreover, $\sigma$ is continuous on $\mathbb{R}\times X$ (see e.g.
\cite{Sch}). The map $\eta$ defined on $X$ by
$\eta(u):=\sigma(2\varepsilon, u)$ satisfies $(i)$. Since
\begin{eqnarray}\label{2.2}
\frac{d}{dt}\varphi\big(\sigma(t,
u)\big)=\bigg(\nabla\varphi\big(\sigma(t, u)\big),
\frac{d}{dt}\sigma(t, u)\bigg)
=\bigg(\nabla\varphi\big(\sigma(t, u)\big), f\big(\sigma(t,
u)\big)\bigg)
=\psi\big(\sigma(t, u)\big).
\end{eqnarray}
If
\begin{eqnarray*}
\sigma(t, u)\in\varphi^{-1}\big([h-\varepsilon,
h-\frac{1}{2}\varepsilon]\big)=B, \quad \ \ \forall \ t\in[0,
2\varepsilon],
\end{eqnarray*}
then
\begin{eqnarray*}
\psi(\sigma(t, u))=1.
\end{eqnarray*}
Let $u\in\varphi^{-1}\big([h-\varepsilon,
h-\frac{1}{2}\varepsilon]\big)$, we obtain from (\ref{2.2}),
\begin{eqnarray*}
\varphi\big(\sigma(2\varepsilon,
u)\big)=\varphi(u)+\int_{0}^{2\varepsilon}\frac{d}{dt}\varphi\big(\sigma(t,
u)\big)dt
=\varphi(u)+\int_{0}^{2\varepsilon}\psi\big(\sigma(t, u)\big)dt
\geq h-\varepsilon+2\varepsilon=h+\varepsilon,
\end{eqnarray*}
and
\begin{eqnarray*}
\varphi\big(\sigma(2\varepsilon,
u)\big)=\varphi(u)+\int_{0}^{2\varepsilon}\frac{d}{dt}\varphi\big(\sigma(t,
u)\big)dt=\varphi(u)+\int_{0}^{2\varepsilon}\psi\big(\sigma(t, u)\big)dt\leq h-\frac{1}{2}\varepsilon+2\varepsilon=h+\frac{3}{2}\varepsilon.
\end{eqnarray*}
So, $(ii)$ is also proved.

Finally, similar to the proof for $(ii)$, we can prove $(iii)$. $\Box$

\begin{remark}\label{2.1}
By using a special domain $D$, a new locally Lipschitz function $\psi$ and some skills, we obtain Lemma \ref{lem2.1}, and in Lemma \ref{lem2.1}, more fixed points, especially for $u=\varphi^{-1}(c)$, are obtained.
\end{remark}

\begin{remark}\label{2.2}
The deformations in Lemma \ref{lem2.1} \big(i.e. conclusion $(ii)$ and $(iii)\big)$ are different from those in Lemma \ref{lem1.1}.
\end{remark}
 \vskip1mm\par
Now, by Lemma \ref{lem2.1}, we can prove our new mountain pass theorem which is independent of $c_0$.

\section{Proof of Theorem \ref{the1.5}}\setcounter{equation}{0}
\vskip1mm\noindent
\textbf{Proof.} Conclusion $(i)$ is obvious. Suppose that conclusion $(ii)$ does not hold. Consider $\beta=\eta\circ\gamma$, where $\eta$ is given by Lemma \ref{lem2.1}, and then we need to check two cases.

\textbf{Case 1.} $\varphi(0)<\varphi(e)=\varphi(e_1)<\hat{c}$.

By an analogue argument of conclusion $(i)$ of Lemma \ref{lem1.1}, we have
 \begin{eqnarray*}
\beta(0)&=&\eta\big(\gamma(0)\big)=\eta(0)=0,\\
\beta(e)&=&\eta\big(\gamma(\frac{1}{2})\big)=\eta(e)=e,\\
\beta(e_1)&=&\eta\big(\gamma(1)\big)=\eta(e_1)=e_1,
     \end{eqnarray*}
So, $\beta\in\hat{\Gamma}$. By the definition of $\hat{c}$, there exists $\gamma\in\hat{\Gamma}$ such that
\begin{eqnarray}\label{3.1}
   \max_{t\in[0, 1]}\varphi\big(\gamma(t)\big)\leq
    \hat{c}+\varepsilon.
 \end{eqnarray}
It follows from conclusion $(ii)$ of Lemma \ref{lem1.1} and (\ref{3.1}) that $\hat{c}\leq\max_{t\in[0, 1]}\varphi\big(\beta(t))\big)\leq \hat{c}-\varepsilon$. This is a contradiction.

\textbf{Case 2.} $\varphi(0)<\varphi(e)=\varphi(e_1)=\hat{c}$.

If $\max_{t\in[0, 1]}\varphi\big(\gamma(t))\big)\equiv \hat{c}$ for $\forall \gamma\in \hat{\Gamma}$, then the theorem is obviously. Combining the definition of $\hat{c}$, there exists $\gamma\in\hat{\Gamma}$ such that
\begin{eqnarray}\label{3.2}
   \hat{c}+\frac{1}{2}\varepsilon\leq\max_{t\in[0, 1]}\varphi\big(\gamma(t))\big)\leq
    \hat{c}+\varepsilon.
 \end{eqnarray}
Take  $D = \{u\in X\mid h=\varphi(u)=\hat{c}\}$ in Lemma \ref{lem2.1}. Then, by conclusion $(i)$ of the new quantitative deformation lemma, we have
 \begin{eqnarray*}
\beta(0)&=&\eta\big(\gamma(0)\big)=\eta(0)=0,\\
\beta(e)&=&\eta\big(\gamma(\frac{1}{2})\big)=\eta(e)=e,\\
\beta(e_1)&=&\eta\big(\gamma(1)\big)=\eta(e_1)=e_1.
     \end{eqnarray*}
So, $\beta\in\hat{\Gamma}$. It follows from the conclusion $(iii)$ of the new quantitative deformation lemma and (\ref{3.2}) that
 $\hat{c}\leq\max_{t\in[0, 1]}\varphi(\beta(t))\leq \hat{c}-\varepsilon$. It is a contradiction.

Combining Case 1 and Case 2, the proof for our new mountain pass theorem is complete. $\Box$

\section{Proof of Theorem \ref{the1.6}}\setcounter{equation}{0}
\vskip1mm\noindent
\textbf{Proof}. We prove Theorem \ref{the1.6} by five steps.

Step 1: We make some notations.
\begin{itemize}
 \item [$\bullet$]  For $a,
\,b  \in \mathbb{Z}$, define $\mathbb{Z} [a] = \{a, a+1, \ldots \}$,
$\mathbb{Z} [a, b] = \{a, a+1, \ldots,
 b\}$ when $a\leq b$.
\item [$\bullet$] Let $S$ be the set of sequences, i.e.
$S=\big\{u=\{u_{n}\}= (_{\cdots}, u_{-n}, _{\cdots}, u_0, _{\cdots},
u_{n}, _{\cdots}), \ u_{n}
 \in \mathbb{R}, \ n \in \mathbb{Z}\big\}$. For any given positive integer $M$, $E_M$
 is defined as a subspace of $S$ by
 \begin{eqnarray*}
E_{M} = \big\{u=\{u_{n}\} \in S \mid u_{n+M} = u_{n}, \ n \in
\mathbb{Z}\big\}.
\end{eqnarray*}
\item [$\bullet$] For $x,y \in S$, $a, b \in \mathbb{R}$, $ax+by$ is defined by
\begin{eqnarray*}
ax+by = \{ax_{n}+by_{n}\}_{n=-\infty}^{+\infty},
\end{eqnarray*}
then $S$ is a vector space. Clearly, $E_{M}$ is isomorphic to
$\mathbb{R}^{M}$, $E_M$ can be equipped with inner product
 \begin{eqnarray*}
 \langle x, y\rangle_{E_M} =
 \sum_{s=1}^{M}x_{s}y_{s},
 \quad \forall \
 x, y \in E_M,
\end{eqnarray*}
then $E_{M}$ with the inner product given above is a finite
dimensional Hilbert space and linearly homeomorphic to
$\mathbb{R}^{M}$. And the norms $\| \cdot \|$ and $\| \cdot
\|_{\beta}$
 induced by
 \begin{eqnarray*}
\|x\| = \bigg(\sum_{j=1}^{M}x_{j}^{2}\bigg)^{\frac{1}{2}},\ \|x\|_{\beta} =
\bigg(\sum_{j=1}^{M}|x_{j}|^{\beta}\bigg)^{\frac{1}{\beta}}, \ \beta\in
[1,\infty),
\end{eqnarray*}
are equivalent, i.e., there exist constants  $0 < C_{1} \leq C_{2}$
such that
\begin{eqnarray*}
C_{1}\| x\| \leq \|x\|_{\beta} \leq C_{2}\|x\|, \quad \forall \ x\in
E_M.
\end{eqnarray*}
\item [$\bullet$] For a given matrix
\begin{eqnarray*}
B = \left(
     \begin{array}{cccccc}
  2 & -1 & 0 & \cdots & 0 & -1 \\
  -1 & 2 & -1 & \cdots & 0 & 0 \\
  0 & -1 & 2 & \cdots & 0 & 0 \\
  \vdots & \vdots & \ddots & \vdots & \vdots & \vdots \\
  -1 & 0 & 0 & \cdots & -1 & 2 \\
 \end{array}
   \right)_{M \times M},
\end{eqnarray*}
then by the results in \cite{Z04}, we have all the eigenvalues of $B$ are 0, $\lambda_1, \,\lambda_2,\,\ldots,\,\lambda_{M-1}$ and $\lambda_j>0$
for all $j\in\mathbb{Z}[1,\,M-1]$. Moreover,
\begin{eqnarray*}
\lambda_{\min} =2(1-\cos\frac{2\pi}{M}), \ \ \ \ \ \ \
 \ \lambda_{\max} =\cases{\ \ \ 4, \ \ \   \mbox{when} \  M \  \mbox{is even}, \cr
2(1+\cos{\frac{\pi}{M}}), \  \mbox{when}  \ M \  \mbox{is odd.}}
\end{eqnarray*}
\end{itemize}

Step 2: Let the functional
\begin{eqnarray}\label{4.1}
 \varphi(u)=
\frac{1}{2}\sum_{s=1}^{M}(\Delta
u_{s})^{2}- F(n, u_{n})-G,
\end{eqnarray}
where
\begin{eqnarray*}
G=G(u_{1}, u_{2}, _{\cdots}, u_{n-1}, u_{n+1}, u_{n+2}, _{\cdots},
u_{M})=w_3\bigg[ \sum_{s=1}^{n-1}
u^{2}_{s}+\sum_{s=n+1}^{M}u^{2}_{s}\bigg].
\end{eqnarray*}
Then, by condition $(A_3)$, we say $\varphi(u)$ is bounded from above on $E_M$. In fact,
according to condition $(A_3)$, if we let $w=\max\{|F(n,\,x)-w_3x^{2}+w_2|: \,n\in\mathbb{Z},\, |x|\leq w_1\}$ and $w^{\prime}=w+w_2$,
then $F(n,\,x)\geq w_3|x|-w^{\prime}$. Combining
$\sum_{s=1}^{M}(\Delta u_{s})^{2}=\sum_{s=1}^{M}(u_{s+1}-u_{s})^{2}=
\sum_{s=1}^{M}(2u_{s}^{2}-2u_{s}u_{s+1})$, we have for all $u\in E_{M}$,
\begin{eqnarray*}
\varphi(u) &=& \frac{1}{2}\bigg[\sum_{s=1}^{M}(\Delta
u_{s})^{2}\bigg]- F(n, u_{n})-G\cr &\leq&
\frac{1}{2}\sum_{s=1}^{M}(2u_{s}^{2}-2u_{s}u_{s+1}) -w_3
u^{2}_n +w^{\prime}-w_3 \sum_{s=1}^{n-1}
u^{2}_{s}-w_3\sum_{s=n+1}^{M}u^{2}_{s}\cr
 &=& \frac{1}{2}u^\top
Bu-w_3\|
u \|^{2} +w^{\prime}
 \cr &\leq&  \frac{1}{2}\lambda_{\max}\|
u \|^{2}-w_3\|
u \|^{2} +w^{\prime}\cr
&=& \big(\frac{1}{2}\lambda_{\max}-w_3\big)\|
u \|^{2} +w^{\prime}.
\end{eqnarray*}
Then, by condition $(A_3)$, we have $\varphi(u)\leq w^{\prime}$. So, $\varphi(u)$ is bounded from above on $E_M$.

Step 3: We claim that $\varphi(u)$ satisfies the $(P.S.)$ condition. In fact, let
$u^{(k)} \in E_{M}$, for all $k\in \mathbb{N}$,
be such that $\{\varphi(u^{(k)})\}$ is bounded. Then, by Step 2, there exists $M_{1} > 0$, such that
\begin{eqnarray*}
-M_{1} \leq \varphi(u^{(k)}) \leq
\big(\frac{1}{2}\lambda_{\max}-w_3\big)\|
u^{(k)} \|^{2} +w^{\prime},
\end{eqnarray*}
which implies that
\begin{eqnarray*}
\|
u^{(k)} \|^{2}\leq(w_3-\frac{1}{2}\lambda_{\max})^{-1}(w_2+M_1).
\end{eqnarray*}
That is, $\{u^{(k)}\}$ is bounded in $E_{M}$. Since $E_{M}$ is
finite dimensional, there exists a subsequence of $\{ u^{(k)} \}$
(not labeled), which is convergent in $E_{M}$, so the $(P.S.)$ condition
is verified.

Step 4: Obviously, from (\ref{4.1}) and condition $(W_2)$, we have $\varphi(0)=0$. Take
\begin{eqnarray*}
e =\cases{ u_{n+1}=\sqrt{w_3}w_4,\,u_{n+2}=-\sqrt{w_3}w_4, \,u_{n-1}=\sqrt{w_3}w_4, \ \ \ \  \  \cr
u_i=0, \    \ i=1,2,\ldots,n-2,n,n+3,\ldots,M,}
\end{eqnarray*}
and
\begin{eqnarray*}
e_1 =\cases{ u_{n+1}=\sqrt{w_3}w_4,\,u_{n+2}=-\sqrt{w_3}w_4,  \ \ \ \  \  \cr
u_i=0, \    \ i=1,2,\ldots,n-1,n,n+3,\ldots,M,}
\end{eqnarray*}
where $w_4>0$.
Then, it is easy to verify that
\begin{eqnarray*}
\varphi(e) &=& \frac{1}{2}\bigg[\sum_{s=1}^{M}(\Delta
u_{s})^{2}\bigg]- F(n, u_{n})-G\cr
 &=& \frac{1}{2}\bigg[(u_{n+3}-u_{n+2})^{2}+(u_{n+2}-u_{n+1})^{2}+(u_{n+1}-u_{n})^{2}+(u_{n}-u_{n-1})^{2}+(u_{n-1}-u_{n-2})^{2}\bigg]\cr&&-0-
 (u_{n+2}^{2}+u_{n+1}^{2}+u_{n-1}^{2})\cr &=&  \frac{1}{2}\bigg[(0+\sqrt{w_3}w_4)^{2}+(-\sqrt{w_3}w_4-\sqrt{w_3}w_4)^{2}+(\sqrt{w_3}w_4-0)^{2}+(0-\sqrt{w_3}w_4)^{2}
 \cr&&+(\sqrt{w_3}w_4-0)^{2}\bigg]-0-
 (w_3w_{4}^{2}+w_3w_{4}^{2}+w_3w_{4}^{2})=w_3w_{4}^{2},
\end{eqnarray*}
and
\begin{eqnarray*}
\varphi(e_1)  &=& \frac{1}{2}\bigg[\sum_{s=1}^{M}(\Delta
u_{s})^{2}\bigg]- F(n, u_{n})-G\cr
 &=& \frac{1}{2}\bigg[(u_{n+3}-u_{n+2})^{2}+(u_{n+2}-u_{n+1})^{2}+(u_{n+1}-u_{n})^{2}\bigg]-0-
 (u_{n+2}^{2}+u_{n+1}^{2})\cr &=&  \frac{1}{2}\bigg[(0+\sqrt{w_3}w_4)^{2}+(-\sqrt{w_3}w_4-\sqrt{w_3}w_4)^{2}+(\sqrt{w_3}w_4-0)^{2}\bigg]-0-
 (w_3w_{4}^{2}+w_3w_{4}^{2})\cr &=&w_3w_{4}^{2}.
\end{eqnarray*}
So, $\varphi(e)=\varphi(e_1)=w_3w_{4}^{2}>0=\varphi(0)$. Moreover, all the conditions of our new mountain pass theorem are satisfied. Noticing that $\varphi(u)$ satisfies the $(P.S.)$ condition, there exists a critical point $\hat{u}$ such that $\varphi(\hat{u})=\hat{c}$ ($\hat{c}$ is given in Theorem \ref{the1.5}).

  Step 5: We say the system (\ref{1.1}) has at least one nontrivial $M$-periodic solutions. In fact,
$\varphi\in C^{2} (E_{M}, \mathbb{R})$. For any $u =
\{u_{n}\}_{n \in\, \mathbb{Z}} \in E_{M}$, according to $u_{0} =
u_{M}$, $u_{1} = u_{M+1}$, one computes that
\begin{eqnarray*}
\frac{\partial \varphi}{\partial
u_{n}}=\Delta^{2}u_{n-1}+\nabla_{u_{n}} F(n, u_{n}),\quad \forall \ n \in \mathbb{Z}.
\end{eqnarray*}
Therefore, the existence of critical points of $\varphi$ on $E_{M}$
implies the existence of periodic solutions of system
(\ref{1.1}).

  Note that when $u_{1}=\ldots=u_{M}$, then $\Delta
u_{1}=\ldots=\Delta u_{M}=0$. Combining (\ref{4.1}) and
$F\geq0$ in condition $(W_1)$, then $\varphi(u)\leq0$. But $\hat{c}\geq\varphi(e)>\varphi(0)=0$, so the above periodic solution
$\hat{u}$ is nontrivial. From this Theorem \ref{the1.6} is proved. $\Box$

\begin{remark}\label{rem4.1}
Take $F(t,x)$ as in Example \ref{exa1.2}, from (\ref{4.1}), we have
\begin{eqnarray*}
\varphi(u)=
\frac{1}{2}\bigg[\sum_{s=1}^{M}(\Delta
u_{s})^{2}\bigg]- a(\mu u^{2}_{n}+\cos u_{n}-1)(\phi(n)+K)-w_3\bigg[ \sum_{s=1}^{n-1}
u^{2}_{s}+\sum_{s=n+1}^{M}u^{2}_{s}\bigg],
\end{eqnarray*}
where $u= (_{\cdots}, u_{-n}, _{\cdots}, u_0, _{\cdots})\in E_{M}$. We notice that the value of $\inf_{\|u\|=r}\varphi(u)$ is very difficult to compute, but fortunately, the condition in our new mountain pass theorem (Theorem \ref{the1.5}) is independent of $\inf_{\|u\|=r}\varphi(u)$, we need not compute it.
\end{remark}

\vskip6mm\noindent
\textbf{\large{Acknowledgements}}
 \vskip4mm\par
We sincerely thank the editors and referees for their
valuable comments.

\end{document}